\documentclass[sn-mathphys-num]{sn-jnl}% Math and Physical Sciences Numbered Reference Style
\usepackage{graphicx}%
\usepackage{multirow}%
\usepackage{amsmath,amssymb,amsfonts}%
\usepackage{amsthm}%
\usepackage{mathrsfs}%
\usepackage[title]{appendix}%
\usepackage{xcolor}%
\usepackage{textcomp}%
\usepackage{manyfoot}%
\usepackage{booktabs}%
\usepackage{algorithm}%
\usepackage{algorithmicx}%
\usepackage{algpseudocode}%
\usepackage{listings}%
\theoremstyle{thmstyleone}%
\newtheorem{theorem}{Theorem}%  meant for continuous numbers
\theoremstyle{thmstyletwo}%
\newtheorem{lemma}{Lemma}

\theoremstyle{thmstylethree}%

\makeatletter
\newcommand{\lastemail}[1]{%
	\global\advance\emailcnt by 1\relax%
	\if@corauemail%
	\g@addto@macro\corrauthemail{%
		\setcounter{footnote}{0}%
		\textcolor{blue}{#1}%
	}%
	\else%
	\g@addto@macro\authemail{%
		\setcounter{footnote}{0}%
		\textcolor{blue}{#1}%
	}%
	\fi
}
\makeatother

\begin{document}

\title[Nevanlinna theory of the Hahn difference operators]{Nevanlinna theory of the Hahn difference operators and its applications}

\author[1]{ \sur{Ling Wang}}\email{wangling812@gznu.edu.cn}

\author*[1]{\sur{Jianren Long}}\email{longjianren2004@163.com}
%\equalcont{These authors contributed equally to this work}
\author[1]{\sur{Xuxu Xiang}}\lastemail{xiangxuxu@gznu.edu.cn}

\affil[1]{\orgdiv{School of Mathematical Sciences}, \orgname{Guizhou Normal University}, \orgaddress{ \city{Guiyang}, \postcode{550025}, \country{China}}}


\abstract{The version of Nevanlinna theory based on Hahn difference operators $$\mathcal{D}_{q,c}(g)=\frac{g(qz+c)-g(z)}{(q-1)z+c}$$ for meromorphic function \(g\) of zero order is studied in this paper. Firstly, the logarithmic derivative lemma for the Hahn difference operators is obtained. As applications of the logarithmic derivative lemma, the Clunie lemma for the Hahn difference operators and the growth of solutions to Hahn difference equations are considered, respectively. Secondly, the second main theorem for the Hahn difference operators is built which is utilized to prove the deficiency relation, Picard's theorem, and the five-value theorem in the sense of Hahn difference operators. Finally, the existence of solutions to Fermat-type Hahn difference equations is also discussed.}



%\pacs
\maketitle

\section{Introduction}
\indent In 1949, Hahn \cite{Hahn} introduced a class of difference operators, which were defined as
\begin{equation*}
	\mathcal{D}_{q,c}g(z)=\frac{g(qz+c)-g(z)}{(q-1)z+c},
\end{equation*}
where $q\in\mathbb{C}\setminus\{0,1\}$, $c$ is a complex constant, \(g\) is a meromorphic function, and the difference operator is called Hahn difference operator, in which the orthogonal polynomials satisfying q-difference equations was first time systematically investigated, thereby laying the foundation for the entire field. Subsequently, Hamza \cite{Ham} investigated Sturm Liouville problems in the setting of Hahn differences, representing a significant step in extending the theory of differential equations to the discrete case. Hamza et al. \cite{Hamza} systematically extended the theory of Hahn difference equations from the scalar case to Banach algebras, establishing the existence and uniqueness of solutions, the definitions and properties of exponential and trigonometric functions, the Wronskian theory, and solution methods for non-homogeneous equations.\\
\indent It is noteworthy that the Hahn difference operator, as a mathematical tool unifying the \(q\)-difference (scale transformation) and the classical difference (translation transformation), has demonstrated significant application value across multiple disciplinary fields in recent years. For instance, the Hahn difference operator has been successfully applied to the modeling of Dirac systems in quantum mechanics. H\i ra \cite{H, Hira} established a sampling theorem for \(q,w\)-Dirac systems associated with the Hahn difference operator, providing a theoretical foundation for signal reconstruction in discrete quantum systems. Research indicates that the Hahn difference operator can accurately describe the evolutionary behavior of quantum systems with non-uniform lattice structures, offering a new mathematical tool for addressing electron states in lattice defects or aperiodic potential fields.

The Hahn difference operator also plays an important role in quantum variational problems. Malinowska and Torres \cite{Mal, M} systematically investigated the Hahn quantum variational principle and established a variational method applicable to discontinuous functions. The core advantage of this method lies in its ability to handle discontinuity problems that traditional variational methods cannot address, which is particularly crucial for describing physical processes with discontinuous characteristics, such as quantum jumps and critical phenomena in phase transitions.

Furthermore, the sampling theorem serves as the cornerstone of digital signal processing, and the Hahn difference operator plays a key role in the reconstruction of non-uniformly sampled signals. Annaby \cite{Anna} established a sampling theorem for the Jackson-N\"orlund transform associated with the Hahn difference operator, which generalizes the classical Shannon sampling theorem to signal spaces based on the Hahn difference operator. Their findings reveal that when a signal satisfies the orthogonal system defined by the Hahn difference operator, accurate reconstruction can be achieved based on non-uniform time nodes. This result holds broad application prospects in areas such as non-uniformly sampled signal processing in communication systems, variable-rate data transmission, and discrete representation of band-limited signals.

Almutairi \cite{Alm} systematically investigated the Hyers-Ulam stability, Rassias stability, and Mittag-Leffler-type stability of nonlinear quantum difference equations based on the \(\beta\)-difference operator. Their research demonstrates that the Hahn difference equation framework can unify the stability analysis of classical difference equations and \(q\)-difference equations, and is particularly suitable for describing discrete dynamical systems with both scale-dependent and translation-dependent characteristics. This theory has significant application value in the stability analysis of digital filter design, adaptive control systems, and networked control systems.\\
\indent These applications demonstrate that the Hahn difference operator exhibits broad application prospects in the modeling and analysis of practical problems. Naturally, this raises the question: can the Hahn difference operator also be applied in the complex domain? To address this question in the future, we must first investigate what properties the Hahn difference operator possesses in the complex domain.\\
\indent Recently, Zhang et al. \cite{Zhang} utilized Nevanlinna theory to investigate the value sharing problem between entire function \(g\) with Picard exceptional values and \(\mathcal{D}_{q,c}g\), which developed the results of the Br\"{u}ck conjecture. However, the Nevanlinna theory for Hahn difference operators has not been systematically studied.\\	
\indent It is worth noting that Hahn difference operators unifies q-Jackson difference operators and forward difference operators, which are defined by
\begin{equation*}
	\mathcal{D}_{q}g(z)=\frac{g(qz)-g(z)}{(q-1)z}
\end{equation*}
and
\begin{equation*}
	\Delta^{c}g(z)=\frac{g(z+c)-g(z)}{c},
\end{equation*}
respectively. Therefore, it is evident that
\begin{align*}
	\lim\limits_{c\rightarrow0}\mathcal{D}_{q,c}g(z)=\mathcal{D}_{q}g(z);&~~~~~\lim\limits_{q\rightarrow1}\mathcal{D}_{q,c}g(z)=\Delta^{c}g(z);\\
	\lim\limits_{c\rightarrow0,q\rightarrow1}\mathcal{D}_{q,c}g(z)=g'(z);& \lim\limits_{c\rightarrow1,q\rightarrow1}\mathcal{D}_{q,c}g(z)=\Delta g(z);
\end{align*}	
\indent It is well-known that the Nevanlinna theory based on the classical differential operators \(g^{(k)}\) was established by Nevanlinna in the 1920s, where \(k\) is any positive integer. It has played the key role in studying the properties of solutions of complex differential equations. Later on, the Nevanlinna theory for some difference operators was investigated. For instance, the Nevanlinna theory for classical shift difference operators $\Delta_{c}(g)=g(z+c)-g(z), c\in\mathbb{C}\setminus\{0\}$ was discussed firstly by Halburd-Korhonen \cite{HK1, HK2} and Chiang-Feng \cite{CF1,CF2}, independently. Besides, Barnett et al. \cite{BH} considered the Nevanlinna theory of the classical q-difference operators $\Delta_{q}(g)=g(qz)-g(z), q\in\mathbb{C}\setminus\{0,1\}$. Recently, Cao et al. \cite{CD} established the Nevanlinna theory for q-Jackson difference operators \(\mathcal{D}_{q}g(z)\). Motivated by the rich theoretical background and recent advances in Nevanlinna theory for different operators, we investigate Nevanlinna theory for the Hahn difference operators \(\mathcal{D}_{q,c}g(z)\) by employing techniques developed for both classical difference operators $\Delta_{c}g(z) = g(z+c)-g(z)$ and $q$-difference operators $\Delta_{q}g(z) = g(qz)-g(z)$.\\
\indent For the convenience of readers, it is necessary to introduce some essential notations. For $a\in\mathbb{C},n\in\mathbb{N}$, then
\begin{align*}
	(a;q)_{0}=1,(a;q)_{n}=(1-a)(1-aq)(1-aq^{2})\cdots(1-aq^{n-1}).
\end{align*}
If $0<|q|<1$, then $(a;q)_{\infty}=\prod_{n=0}^{\infty}(1-aq^{n})$. Furthermore,
\[
\begin{pmatrix}
	n \\
	j
\end{pmatrix}_{q}=\frac{(q;q)_{n}}{(q;q)_{j}(q;q)_{n-j}}.
\]
\indent Assuming the reader is familiar with the basic notation and main results of Nevanlinna theory. For example, proximity function $m(r,g)$, counting function $N(r,g=\alpha), \alpha\in\mathbb{C}\cup\{\infty\}$ and characteristic function $T(r,g)$. \\
\indent For a meromorphic function \(g\) in the complex plane \(\mathbb{C}\), the order of  \(g\) is defined as
\begin{align*}
	\rho(g)=\limsup\limits_{r\to\infty}\frac{\log^{+}T(r,g)}{\log r}.
\end{align*}
\indent This paper is organized as follows. In Section 2, the logarithmic derivative lemma and its applications are presented, such as Clunie lemma and the growth of solutions of Hahn difference equations. The second main theorem and its applications are studied in Section 3, such as Defect relation, Picard theorem and Five-values theorem. The existence of solutions for Fermat-type Hahn difference equations is investigated in Section 4. The final part proposes an open question for further discussion.
\section{Logarithmic derivative lemma for Hahn difference operators and its application}

\indent In this section, we will consider two versions of the logarithmic derivative lemma for the Hahn difference operators, based on the difference logarithmic derivative lemma and $q$-difference logarithmic derivative lemma, and denote $k$-th order Hahn difference operators as
\begin{align*}
	\mathcal{D}_{q,c}^{k}g(z)=\mathcal{D}_{q,c}(\mathcal{D}_{q,c}^{k-1}g(z)); \mathcal{D}_{q,c}^{0}g(z)=g(z), k\geq1.
\end{align*}

%\subsection{Logarithmic derivative lemma I}
\indent In order to prove Theorem \ref{thm3.1} in the following, we should introduce some lemmas. The first lemma is a summary of [4, Theorem 1.1] and [47, Theorem 1.1].
\begin{lemma}\cite{BH, zhang}\label{l3.1}
	Let $g$ be a non-constant meromorphic function with zero order, and $q\in\mathbb{C}\setminus\{0\}$. Then for all $r\in E_{1}$,
	\begin{align*}\label{l3.1'}
		m\left(r,\frac{g(qz)}{g(z)}\right)=o(T(r,g)),\quad T(r,g(qz))=T(r,g)+o(T(r,g)),
	\end{align*}
	where the set $E_{1}$ is of logarithmic density 1.
\end{lemma}
\indent The second lemma is the logarithmic derivative lemma for \(g(z+c)\), which can be found in [14, Corollary 1.2.1].
\begin{lemma}\cite{HK1}\label{l3.2}
	Let $g$ be a non-constant meromorphic function with finite order, let $c\in\mathbb{C}$ and $\delta<1$. Then
	\begin{align*}%\label{l3.2'}
		m\left(r,\frac{g(z+c)}{g(z)}\right)=o\left(\frac{T(r+|c|,g)}{r^{\delta}}\right)=o\left(\frac{T(r,g)}{r^{\delta}}\right)=o(T(r,g))
	\end{align*}
	holds for all $r\not\in E_{2}$, where the set $E_{2}$ is of finite logarithmic measure.
\end{lemma}
\indent Similarly, we also need to demonstrate the relationship between the characteristic functions of \(g(z+c)\) and \(g\). Hence, we present the following conclusion, which originates from [8, Theorem 2.1].
\begin{lemma}\cite{CF1}\label{l3.3}
	Let $g$ be a meromorphic function with finite order \(\rho(g)\), and let \(c\) be a fixed non-zero complex number. Then for each \(\varepsilon>0\), we have
	\begin{align*}%\label{l3.3'}
		T(r,g(z+c))=T(r,g)+O(r^{\rho(g)-1+\varepsilon})+O(\log r).
	\end{align*}
\end{lemma}

\indent Next, we obtain the first result as below.

\begin{theorem}\label{thm3.1}
	Suppose that $g$ is a non-constant meromorphic function with zero order. Then for any $0<|q|<1$,
	\begin{align*}
		m\left(r,\frac{\mathcal{D}_{q,c}g(z)}{g(z)}\right)=o(T(r,g)),
		m\left(r,\frac{\mathcal{D}_{q,c}^{k}g(z)}{g(z)}\right)=o(T(r,g))
	\end{align*}
	hold on a set of logarithmic density 1.
\end{theorem}

\begin{proof}
	Applying Lemmas \ref{l3.1}-\ref{l3.3}, then there exists a set $E=E_{1}\setminus E_{2}$ with logarithmic density 1, such that for all $r\in E$,
	\begin{align*}%\label{t1}
		m\left(r,\frac{\mathcal{D}_{q,c}g(z)}{g(z)}\right)&\leq m\left(r,\frac{g(qz+c)}{g(z)}\right)+o(T(r,g))\nonumber\\
		&\leq m\left(r,\frac{g(qz+c)}{g(qz)}\right)+m\left(r,\frac{g(qz)}{g(z)}\right)+o(T(r,g))\\
		&=o(T(r,g)).\nonumber
	\end{align*}
	
	For any positive integer $k\geq2$, it follows from the definition of $\mathcal{D}_{q,c}^{k}g(z)$ that
	\begin{equation}\label{m1}
		\mathcal{D}_{q,c}^{k}g(z)=\frac{1}{((q-1)z+c)^{k}}\sum_{i=0}^{k}(-1)^{i}
		\begin{pmatrix}
			n \\
			j
		\end{pmatrix}_{q}q^{\frac{i(i-1)}{2}}g\left(q^{k-i}z+c\cdot\frac{q^{k-i}-1}{q-1}\right).
	\end{equation}
	\indent Then,	
	\begin{align}\label{m2}
		&~~~~~~~~~~~m\left(r,\frac{g(q^{k-i}z+c\cdot\frac{q^{k-i}-1}{q-1})}{g(z)}\right)\nonumber\\
		&~~~~~~~~~~~\leq\sum_{n=0}^{k-i-1}m\left(r,\frac{g(q^{n+1}z+c\cdot\frac{q^{n+1}-1}{q-1})}{g(q^{n}z+c\cdot\frac{q^{n}-1}{q-1})}\right)\\
		 &\leq\sum_{n=0}^{k-i-1}\left(m\left(r,\frac{g(q^{n+1}z+\eta_{1})}{g(q^{n+1}z)}\right)+m\left(r,\frac{g(q^{n+1}z)}{g(q^{n}z)}\right)+m\left(r,\frac{g(q^{n}z)}{g(q^{n}z+\eta_{2})}\right)\right)\nonumber\\
		&=o(T(r,g)),\nonumber
	\end{align}
	where $\eta_{1}=c\cdot\frac{q^{n+1}-1}{q-1}, \eta_{2}=c\cdot\frac{q^{n}-1}{q-1}$. \\
	\indent Hence,
	\begin{align*}
		m\left(r,\frac{\mathcal{D}_{q,c}^{k}g(z)}{g(z)}\right)=o(T(r,g))
	\end{align*}
	holds for all \(|z|=r\in E\). Then, Theorem \ref{thm3.1} is proved completely.
\end{proof}
\indent The Clunie lemma \cite{hayman} is often used to study the properties of solutions of non-linear differential equations (e.g., those involving exponential or trigonometric terms), such as the existence of non-constant meromorphic solutions, or the specific forms (such as rational functions, exponential functions, etc.) of solutions. Additionally, the Clunie lemma is one of the key tools for analyzing the growth and value distribution of transcendental solutions of Painlev\'{e} equations. Besides, it is worth noting that there exists a corresponding difference version of the Clunie lemma \cite{chen}. Therefore, we consider the Hahn Clunie lemma as an application of Theorem \ref{thm3.1}.\\
\indent Before stating Theorem \ref{thm2.5} as below, it is necessary to first explain that Hahn difference polynomial \(P(z,g)\) is expressed as
\begin{align*}
	P(z,g)=\sum_{I} a_{I}(z)g^{l_{0}}(\mathcal{D}_{q,c}g(z))^{l_{1}}\cdots(\mathcal{D}^{v}_{q,c}g(z))^{l_{v}},~l_{0}+l_{2}+\cdots+l_{v}=I,	
\end{align*}
where \(a_{I}\) is the small function of \(g\), \(l_{0},...,l_{v}\) are non-negative integers, \(v\) is a finite non-negative integer, and the maximum value of \(I\) is called the degree of the Hahn difference polynomial \(P(z,g)\).
\begin{theorem}\label{thm2.5}
	Suppose that \(g\) is a transcendental meromorphic function with zero order and that
	\begin{align*}
		g^{n}(z)P(z,g)=Q(z,g),
	\end{align*}
	where \(P(z,g), Q(z,g)\) are Hahn difference polynomials in \(g\) and its any orders Hahn difference operators, and the degree of \(Q(z,g)\) is at most \(n\). Then
	\begin{align*}
		m(r,P(z,g))=o(T(r,g))
	\end{align*}
	as \(r\to\infty\).
\end{theorem}
\begin{proof}
	\indent We abbreviate \(P(z,g)\) as \(P(z)\) in our proof. According to the definition of proximity function \(m(r,g)\), then
	\begin{align*}
		m(r,P(z))&=\frac{1}{2\pi}\int_{0}^{2\pi}\log^{+}|P(re^{i\theta})|d\theta\\
		&=\frac{1}{2\pi}\int_{E}\log^{+}|P(re^{i\theta})|d\theta+\frac{1}{2\pi}\int_{E^{c}}\log^{+}|P(re^{i\theta})|d\theta,
	\end{align*}
	where set \(E=\{\theta\in[0,2\pi):|g(re^{i\theta})|<1\}\), and \(E^{c}\) is the complementary set of \(E\). On the set \(E\), denote
	\begin{align*}
		G(z)=a_{I}(z)g(z)^{l_{0}}(\mathcal{D}_{q,c}g(z))^{l_{1}}\cdots(\mathcal{D}^{v}_{q,c}g(z))^{l_{v}},
	\end{align*}
	then,
	\begin{align*}
		|G(z)|\leq|a_{I}|\left|\frac{\mathcal{D}_{q,c}g(z)}{g(z)}\right|^{l_{1}}\cdots\left|\frac{\mathcal{D}^{v}_{q,c}g(z)}{g(z)}\right|^{l_{v}}.
	\end{align*}
	Hence,
	\begin{align*}
		\frac{1}{2\pi}\int_{E}\log^{+}|G(re^{i\theta})|d\theta&\leq m(r,a_{I})+O\left(\sum_{t=1}^{v}m\left(r,\frac{\mathcal{D}^{v}_{q,c}g(z)}{g(z)}\right)\right)\\
		&+o(T(r,g)).
	\end{align*}
	Combining the Theorem \ref{thm3.1}, then
	\begin{align*}
		\frac{1}{2\pi}\int_{E}\log^{+}|P(re^{i\theta})|d\theta=\sum_{I}\int_{E}\log^{+}|G(re^{i\theta})|d\theta+O(1)=o(T(r,g)).
	\end{align*}
	\indent On the other hand, on the set \(E^{c}=\{\theta\in[0,2\pi):|g(re^{i\theta})|\geq1\}\). The Hahn difference polynomial \(Q(z,g)\) is the sum of terms of the form
	\begin{align*}
		a_{I}(z)g(z)^{l_{0}}(\mathcal{D}_{q,c}g(z))^{l_{1}}\cdots(\mathcal{D}^{v}_{q,c}g(z))^{l_{v}},~l_{0}+l_{1}+\cdots+l_{v}=I\leq n,
	\end{align*}
	\(v\) is a non-negative integer. So,
	\begin{align*}
		|P(z)|&=\left|\frac{1}{g^{n}(z)}\sum_{I}a_{I}(z)g(z)^{l_{0}}(\mathcal{D}_{q,c}g(z))^{l_{1}}\cdots(\mathcal{D}^{v}_{q,c}g(z))^{l_{v}}\right|\\
		&\leq\sum_{I}\left|a_{I}(z)\left|\frac{\mathcal{D}_{q,c}g(z)}{g(z)}\right|^{l_{1}}\cdots
		\left|\frac{\mathcal{D}^{v}_{q,c}g(z)}{g(z)}\right|^{l_{v}}\right|.
	\end{align*}
	Hence,
	\begin{align*}
		 \frac{1}{2\pi}\int_{E^{c}}\log^{+}|P(re^{i\theta})|d\theta&\leq\left(\sum_{t=1}^{v}m\left(r,\frac{\mathcal{D}^{t}_{q,c}g(z)}{g(z)}\right)+m(r,a(z))\right)\\
		&=o(T(r,g)).
	\end{align*}
	Therefore, Theorem \ref{thm2.5} is completely proved.
\end{proof}
%\section{Logarithmic derivative lemma for Hahn difference operator and its application II}
\indent For a meromorphic function \(g\) in complex plane \(\mathbb{C}\), its logarithmic order is defined as
\begin{align*}
	\rho_{\log}(g)=\limsup\limits_{r \to \infty}	\frac{\log^{+}T(r,g)}{\log \log r}.
\end{align*}
\indent Next, the another logarithmic derivative lemma for Hahn difference operators is described by using the logarithmic order.
\begin{theorem}\label{l5.1}
	Suppose that \(g\) is a non-constant meromorphic function with finite logarithmic order, \(q'\in\mathbb{C}\setminus\{0\}, c\in\mathbb{C}\), \(k\) is a positive integer. Then for any given \(\varepsilon>0\),
	\begin{align*}
		m\left(r,\frac{\mathcal{D}_{q',c}^{k}g(z)}{g(z)}\right)=O((\log r)^{\rho_{\log}(g)-1+\varepsilon}).
	\end{align*}
\end{theorem}
\begin{proof}
	The logarithmic order of meromorphic function \(g\) is finite, then the order of \(g\) is zero. According to [40, Theorem 2.2] and [8, Corollary 2.6], then
	\begin{align*}
		m\left(r,\frac{g(q'z)}{g(z)}\right)&=O((\log r)^{\rho_{\log}(g)-1+\varepsilon});\\
		m\left(r,\frac{g(q'z+c)}{g(q'z)}\right)&=O(r^{\rho(g)-1+\varepsilon})=O(\frac{1}{r}).
	\end{align*}	
	\indent Considering \eqref{m1} and \eqref{m2}, then
	\begin{align*}
		m\left(r,\frac{\mathcal{D}_{q',c}^{k}g(z)}{g(z)}\right)= O((\log r)^{\rho_{\log}(g)-1+\varepsilon}).
	\end{align*}
\end{proof}

\indent It well known that the growth or other properties of solutions of differential equations and difference equations is very interesting topic, many results have been obtained, such as \cite{CF1, GGG, Hei, Hel, li, l2016, l2018, l2024, l2025, wj2008, wpc2011}. However, the slowing growth solutions of differential equations and difference equations is difficult to study. Here, we consider the slowing growth solutions of the linear Hahn difference equations
\begin{equation}\label{5.1}
	\mathcal{D}_{q,c}^{k}g(z)+A_{k-1}(z)\mathcal{D}_{q,c}^{k-1}g(z)+\cdots+A_{1}(z)\mathcal{D}_{q,c}g(z)+A_{0}(z)g(z)=0,
\end{equation}
where \(|q|\neq0,1, c\in\mathbb{C}, k\) is a positive integer, and \(A_{i}(i=0,...,k-1)\) are zero order entire functions.
\begin{theorem}
	Let \(A_{i}\) be entire functions such that \(\rho_{\log}(A_{0})>\rho_{\log}(A_{j})\geq1,i=0,1,...,k-1,j=1,2,...,k-1\), \(k\geq2\) is an integer. Suppose that \(g\) is an entire solution of equation \eqref{5.1}. Then \(g\) satisfies \(\rho_{\log}(g)\geq\rho_{\log}(A_{0})+1\).
\end{theorem}
\begin{proof}
	\indent By equation \eqref{5.1}, then
	\begin{align}
		T(r,A_{0}(z))&=m(r,A_{0}(z))\nonumber\\
		&\leq T(r,A_{k-1}(z))+\cdots+T(r,A_{1}(z))\\
		&+m\left(r,\frac{\mathcal{D}_{q,c}^{k-1}g(z)}{g(z)}\right)+\cdots+m\left(r,\frac{\mathcal{D}_{q,c}g(z)}{g(z)}\right)\nonumber.
	\end{align}
	\indent Combining Theorem \ref{l5.1} and the definition of \(\rho_{\log}(A_{0})\), then there exists an infinite sequence \(\{r_{n}\}\) satisfying \(\lim\limits_{r\to\infty}r_{n}=\infty\), such that for any given \(\varepsilon\in(0,\frac{\rho_{\log}(A_{0})-\rho_{\log}}{3})\) and sufficiently large \(n\),
	\begin{align}\label{pf5.1}
		(\log r_{n})^{\rho_{\log}(A_{0})-\varepsilon}\leq O((\log r_{n})^{\rho_{\log}+\varepsilon})+O((\log r_{n})^{\rho_{\log}(g)-1+\varepsilon}),
	\end{align}
	where \(\rho_{\log}=\max\{\rho_{\log}(A_{1}),...,\rho_{\log}(A_{k-1})\}\). It can be deduce from \eqref{pf5.1} that
	$$(1-o(1))(\log r_{n})^{\rho_{\log}(A_{0})-\varepsilon}\leq O((\log r_{n})^{\rho_{\log}(g)-1+\varepsilon}),$$
	which implies \(\rho_{\log}(g)\geq\rho_{\log}(A_{0})+1\).
\end{proof}

\section{The second main theorem for Hahn difference operators and its applications}
%\subsection{Hahn difference analogue of second main theorem}
\indent For any $a\in\mathbb{C}\cup\{\infty\}$, $\overline{n}(r,g=a)=\overline{n}(r,\frac{1}{g-a})$ denotes the number of zeros of \(g(z)-a\) in \(\{z:|z|\leq r\}\), each zero is counted only once. In fact, it can also be denoted as a sum of integers $n-m$, where \(n\) is the number of the zeros of $g-a$ in $\{z:|z|\leq r\}$ with multiplicity, and $m(=n-1)$ is the multiplicity of $g'=0$ and $g=a$. Besides,  \(\overline{n}(r,g)=\overline{n}(r,g=\infty)=\overline{n}\left(r,\frac{1}{g}=0\right)\) denotes the number of poles of \(g\) in \(\{z:|z|\leq r\}\) (or the number of zeros of \(\frac{1}{g}\) in \(\{z:|z|\leq r\}\)), each pole of \(g\) (zero of \(\frac{1}{g}\)) is counted only once. Similarly, it can also be written as a sum of integers \(n-m\), where \(n\) denotes the number of poles of \(g\) in \(\{z:|z|\leq r\}\) with multiplicity, and \(m(=n-1)\) is the multiplicity of \(d\left(\frac{1}{g}\right)\textfractionsolidus dz=-g'(z)\textfractionsolidus g^{2}=0\) where \(g=\infty\).

\indent Next, we define a Hahn difference analogue of the $\overline{n}(r,\frac{1}{g-a})$ of $g$ as
\begin{align*}
	\hat{n}_{q,c}(r,g=a)=\hat{n}_{q,c}\left(r,\frac{1}{g-a}\right)
\end{align*}
to be the sum of integers of the form $n-m$, where \(n\) denotes the multiplicity of the point $z$, which satisfies $g=a$ in $\{z:|z|\leq r\}$. Denote $m$ is defined by $m=\min\{n,m'\}$, $m'$ is the multiplicity of the point \(z\), which satisfies $\mathcal{D}_{q,c}g=0$ and $g=a$.

Similarly, we define
\begin{align*}
	\hat{n}_{q,c}(r,g)=\hat{n}_{q,c}(r,g=\infty)=\hat{n}_{q,c}\left(r,\frac{1}{g}=0\right)
\end{align*}
to be the sum of integers of the form $n-m$, where \(n\) denotes the multiplicity of the point $z$, which satisfies $g=\infty$ in $\{z:|z|\leq r\}$. Denote $m$ is defined by $m=\min\{n,m'\}$, $m'$ is the multiplicity of point \(z\), which satisfies $\mathcal{D}_{q,c}\left(\frac{1}{g}\right)=-\frac{\mathcal{D}_{q,c}g}{g(qz+c)g}=0$ and \(g=\infty\).

\indent We also need to define the Hahn~type~integrated counting function of $g$ by
\begin{align*}
	&\hat{N}_{q,c}(r,g=a)=\hat{N}_{q,c}\left(r,\frac{1}{g-a}\right)\\
	&=\int_{0}^{r}\frac{\hat{n}_{q,c}(t,g=a)-\hat{n}_{q,c}(0,g=a)}{t}dt+\hat{n}_{q,c}(0,g=a)\log r
\end{align*}
and
\begin{align*}
	\hat{N}_{q,c}(r,g)=\int_{0}^{r}\frac{\hat{n}_{q,c}(t,g)-\hat{n}_{q,c}(0,g)}{t}dt+\hat{n}_{q,c}(0,g)\log r.
\end{align*}
\indent It is evident that $\hat{N}_{q,c}(r,g=a)$ and $\hat{N}_{q,c}(r,g)$ are counterpart of the $\overline{N}(r,g=a)$ and $\overline{N}(r,g)$ from the classical Nevanlinna theory, respectively, where $\overline{N}(r,g=a)$ is the reduced \(a\)-point counting function, and $\overline{N}(r,g)$ is the reduced pole counting function.\\
\indent The second fundamental theory for Hahn difference operator is established by using Hahn type counting functions defined above.
\begin{theorem}\label{thm2}
	Suppose that $g$ is a non-constant meromorphic function of zero order, $0<|q|<1$, and $a_{1},a_{2},\ldots,a_{l}$ are distinct points in $\mathbb{C}\cup\{\infty\}$, \(l\ge 3\) is integer. Then,
	\begin{align}\label{t1}
		(l-2)T(r,g)&\leq\sum_{i=1}^{l}N(r,g=a_{i})-N_{q,c}(r)+o(T(r,g))\nonumber\\
		&\leq\sum_{i=1}^{l}\hat{N}_{q,c}(r,g=a_{i})+o(T(r,g))
	\end{align}
	holds for all $r\in E_{1}$, where the set $E_{1}$ is of logarithmic density 1, and $$N_{q,c}(r)=2N(r,g)-N(r,\mathcal{D}_{q,c}g)+N\left(r,\frac{1}{\mathcal{D}_{q,c}g}\right).$$
\end{theorem}
\begin{proof}
	Denote
	\begin{align*}
		G(z)=\sum_{v=1}^{l}\frac{1}{g(z)-a_{v}},
	\end{align*}	
	then
	\begin{align}\label{t2.1}
		m(r,G(z))&=m\left(r,G\cdot\mathcal{D}_{q,c}g\cdot\frac{1}{\mathcal{D}_{q,c}g}\right)\nonumber\\
		&\leq m\left(r,\frac{1}{\mathcal{D}_{q,c}g}\right)+m\left(r,\sum_{v=1}^{l}\frac{\mathcal{D}_{q,c}g}{g-a_{v}}\right).
	\end{align}	
	
	On the other hand, for any $v\in\{1,2,\ldots,l\}$, we can rewrite $G(z)$ as follows,
	\begin{align*}
		G(z)=\frac{1}{g(z)-a_{v}}\left(1+\sum_{v=1,\mu\neq v}^{l}\frac{g(z)-a_{v}}{g(z)-a_{\mu}}\right).
	\end{align*}
	
	Let $\delta=\min\{|a_{s}-a_{t}|,s\neq t\}$, and by using the similar calculation in [17, P.32], we can deduce that
	\begin{align}\label{t2.2}
		m(r,G(z))\geq\sum_{v=1}^{l}m(r,a_{v})-l\log\frac{2l}{\delta}-\log2.
	\end{align}
	Combining \eqref{t2.1} and \eqref{t2.2}, then
	\begin{align}\label{t2.3}
		 m\left(r,\frac{1}{\mathcal{D}_{q,c}g}\right)\geq\sum_{v=1}^{l}m(r,a_{v})-m\left(r,\sum_{v=1}^{l}\frac{\mathcal{D}_{q,c}g}{g-a_{v}}\right)-l\log\frac{2l}{\delta}-\log2.
	\end{align}
	
	Adding $N\left(r,\frac{1}{\mathcal{D}_{q,c}g}\right)$ to both sides of \eqref{t2.3}, and applying the first main theorem,
	\begin{align}\label{t2.4}
		T(r,\mathcal{D}_{q,c}g)&=T\left(r,\frac{1}{\mathcal{D}_{q,c}g}\right)+O(1)\nonumber\\
		&=m\left(r,\frac{1}{\mathcal{D}_{q,c}g}\right)+N\left(r,\frac{1}{\mathcal{D}_{q,c}g}\right)+O(1)\\
		&\geq N\left(r,\frac{1}{\mathcal{D}_{q,c}g}\right)+\sum_{v=1}^{l}m(r,a_{v})-m\left(r,\sum_{v=1}^{l}\frac{\mathcal{D}_{q,c}g}{g-a_{v}}\right)+O(1).\nonumber
	\end{align}
	
	In fact,
	\begin{align}\label{t2.5}
		T(r,\mathcal{D}_{q,c}g)&=m(r,\mathcal{D}_{q,c}g)+N(r,\mathcal{D}_{q,c}g)\nonumber\\
		&\leq N(r,\mathcal{D}_{q,c}g)+m(r,g)+m\left(r,\frac{\mathcal{D}_{q,c}g}{g}\right).
	\end{align}
	
	Compare \eqref{t2.4} with \eqref{t2.5} and rearranging the term, then
	\begin{align}\label{t2.6}
		\sum_{v=1}^{l}m(r,a_{v})&\leq m(r,g)+m\left(r,\sum_{v=1}^{l}\frac{\mathcal{D}_{q,c}g}{g-a_{v}}\right)+N(r,\mathcal{D}_{q,c}g)\nonumber\\
		&+m\left(r,\frac{\mathcal{D}_{q,c}g}{g}\right)-N\left(r,\frac{1}{\mathcal{D}_{q,c}g}\right)+O(1).
	\end{align}
	\indent Adding $m(r,g)$ to the both sides of \eqref{t2.6},
	\begin{align}\label{t2.7}
		m(r,g)+\sum_{v=1}^{l}m(r,a_{v})&\leq 2T(r,g)+m\left(r,\sum_{v=1}^{l}\frac{\mathcal{D}_{q,c}g}{g-a_{v}}\right)
		+m\left(r,\frac{\mathcal{D}_{q,c}g}{g}\right)\nonumber\\
		&-\left(2N(r,g)-N(r,\mathcal{D}_{q,c}g)+N\left(r,\frac{1}{\mathcal{D}_{q,c}g}\right)\right)+O(1).
	\end{align}
	By using Theorem \ref{thm3.1}, then there exists a set \(E_{1}\) with logarithmic density 1, for all \(r\in E_{1}\),
	\begin{align*}
		m\left(r,\sum_{v=1}^{l}\frac{\mathcal{D}_{q,c}g}{g-a_{v}}\right)\leq\sum_{v=1}^{l}m\left(r,\frac{\mathcal{D}_{q,c}(g-a_{v})}{g-a_{v}}\right)
		=o(T(r,g)).
	\end{align*}
	Next, adding $N(r,g)+\sum_{v=1}^{l}N(r,a_{v})$ to the both sides of \eqref{t2.7}, then
	\begin{align*}
		(l+1)T(r,g)&\leq 2T(r,g)+N(r,g)+\sum_{v=1}^{l}N(r,a_{v})+m\left(r,\frac{\mathcal{D}_{q,c}g}{g}\right)&\\
		&-\left(2N(r,g)-N(r,\mathcal{D}_{q,c}g)+N\left(r,\frac{1}{\mathcal{D}_{q,c}g}\right)\right)+o(T(r,g))
	\end{align*}
	holds for all \(r\in E_{1}\). That is, for \(r\in E_{1}\), then
	\begin{align}\label{t2.7'}
		(l-2)T(r,g)\leq&\sum_{v=1}^{l}N(r,a_{v})-N_{q,c}(r)+o(T(r,g)),
	\end{align}
	where $N_{q,c}(r)=2N(r,g)-N(r,\mathcal{D}_{q,c}g)+N\left(r,\frac{1}{\mathcal{D}_{q,c}g}\right)$.\\
	\indent For any $a \in\mathbb{C}\cup\{\infty\}$, by the definition of \(\hat{n}_{q,c}(r,g=a)\), it is evident that the difference between \(n(r,g=a)\) and \(\hat{n}_{q,c}(r,g=a)\) will occur at every point \(z_{0}\), which is the zero of \(g-a\) and \(\mathcal{D}_{q,c}g(z)\) in the disc $\{z:|z|\leq r\}$. That is,
	$$n(r,g=a)-\hat{n}_{q,c}(r,g=a)=m\leq m'$$ holds at \(z_{0}\). Similarly, if $a=\infty$, $n(r,g=\infty)-\hat{n}_{q,c}(r,g=\infty)$ at most the number of zeros of $\mathcal{D}_{q,c}\left(\frac{1}{g}\right)$ at which $g$ has a pole in the disc $\{z:|z|\leq r\}$, with due count of multiplicity. In fact,
	\begin{align*}
		\mathcal{D}_{q,c}\left(\frac{1}{g}\right)=-\frac{\mathcal{D}_{q,c}g}{g(qz+c)g},
	\end{align*}
	then its zeros derive from the zeros of $\mathcal{D}_{q,c}g$ or the poles of $g(qz+c)$ and $g$, and it is worth noting that the poles of $\mathcal{D}_{q,c}g$ must be among the poles of $g(qz+c)$ and $g$. Therefore, the number of zeros of $\mathcal{D}_{q,c}\left(\frac{1}{g}\right)$ is no more than the sum of \(C\) and the difference between \(A\) and \(B\), where \(A\) is the sum of number of the poles of \(g(qz+c)\) and \(g\), \(B\) is the number of the poles of $\mathcal{D}_{q,c}g$, and \(C\) is the number of zeros of \(\mathcal{D}_{q,c}g\). Besides, if $z=-\frac{c}{q-1}$ is a pole of $\mathcal{D}_{q,c}g$, we should add 1 to \(|A-B|+C\).\\
	\indent Hence, we can deduce from discussion above that for $l$ distinct value $a_{1},a_{2},\ldots,a_{l}\in\mathbb{C}\cup\{\infty\}$,
	\begin{align}\label{t2.8}
		&\sum_{i=1}^{l}\left(N(r,g=a_{i})-\hat{N}_{q,c}(r,g=a_{i})\right)\nonumber\\
		&\leq N(r,f(z))+N(r,f(qz+c))\\
		&+N\left(r,\frac{1}{\mathcal{D}_{q,c}g(z)}\right)-N(r,\mathcal{D}_{q,c}g(z))+1.\nonumber
	\end{align}
	%By the Lemma \ref{l3.1} and Lemma \ref{l3.3}, there exists a set $E_{1}$ of logarithmic density 1, such that for all $r\in E_{1}$, \eqref{l3.1'} and \eqref{l3.3'} hold.
	Then,
	\begin{align}\label{t2.9}
		\sum_{i=1}^{l}N(r,g=a_{i})&\leq\sum_{i=1}^{l}\hat{N}_{q,c}(r,g=a_{i})+(2+o(1))N(r,f(z))\nonumber\\
		&+N\left(r,\frac{1}{\mathcal{D}_{q,c}g(z)}\right)-N(r,\mathcal{D}_{q,c}g(z))+1\\
		&=\sum_{i=1}^{l}\hat{N}_{q,c}(r,g=a_{i})+N_{q,c}(r)+1.\nonumber
	\end{align}			
	Substitute \eqref{t2.9} into \eqref{t2.7'},
	\begin{align*}
		(l-2)T(r,g)\leq\sum_{i=1}^{l}\hat{N}_{q,c}(r,g=a_{i})+o(T(r,g))
	\end{align*}
	holds for all \(r\in E_{1}\).
\end{proof}
%\section{Applications of the Second main Theorem}
%\subsection{Defect relation for Hahn difference operators}
\indent For any meromorphic function $g$ and complex constant $a\in\mathbb{C}\cup\{\infty\}$, we denote the Nevanlinna defect, ramification index and multiplicity index of $g$ by $\delta(a,g), \Theta(a,g)$ and $\theta(a,g)$, respectively, their definitions are showed as follows,
\begin{align*}
	\delta(a,g)=1-\limsup\limits_{r\rightarrow\infty}\frac{N\left(r,\frac{1}{g-a}\right)}{T(r,g)},
	\Theta(a,g)=1-\limsup\limits_{r\rightarrow\infty}\frac{\overline{N}\left(r,\frac{1}{g-a}\right)}{T(r,g)}
\end{align*}
and
\begin{align*}
	\theta(a,g)=\liminf\limits_{r\rightarrow\infty}\frac{N\left(r,\frac{1}{g-a}\right)-\overline{N}\left(r,\frac{1}{g-a}\right)}{T(r,g)}.
\end{align*}
It is obvious from the classical Nevanlinna's second main theorem that
\begin{align*}
	\sum\limits_{a\in\mathbb{C}\cup\{\infty\}}(\delta(a,g)+\theta(a,g))\leq\sum\limits_{a\in\mathbb{C}\cup\{\infty\}}\Theta(a,g)\leq2.
\end{align*}		
We aim to establish the Hahn difference counterpart of Nevanlinna defect relation, so we should first define the corresponding symbols. For any given meromorphic function $g$ and any complex constant $a\in\mathbb{C}\cup\{\infty\}$. Denote the Hahn-type multiplicity index and the Hahn-type ramification index by $\theta_{q,c}(a,g)$ and $\Theta_{q,c}(a,g)$, respectively, their definitions are as follows,
\begin{align*}
	\theta_{q,c}(a,g)&=\liminf\limits_{r\rightarrow\infty}\frac{N(r,g=a)-\hat{N}_{q,c}(r,g=a)}{T(r,g)},\\
	\Theta_{q,c}(a,g)&=1-\limsup\limits_{r\rightarrow\infty}\frac{\hat{N}_{q,c}(r,g=a)}{T(r,g)}.
\end{align*}		
\indent Base on the Hahn difference analogue of second main theorem, we get the following the Hahn difference defect relation.
\begin{theorem}
	Suppose that $g$ is a non-constant meromorphic function with zero order, and $0<|q|<1$. Then
	 $$\sum\limits_{a\in\mathbb{C}\cup\{\infty\}}(\delta(a,g)+\theta_{q,c}(a,g))\leq\sum\limits_{a\in\mathbb{C}\cup\{\infty\}}\Theta_{q,c}(a,g)\leq2.$$
\end{theorem}		
\begin{proof}
	By the Theorem \ref{thm2}, then there exists a set \(E_{1}\) with logarithmic density 1, such that for all \(r\in E_{1}\),
	\begin{align}\label{33}
		(l-2)T(r,g)\leq\sum_{i=1}^{l}\hat{N}_{q,c}(r,g=a_{i})+o(T(r,g)).
	\end{align}
	\indent Dividing the characteristic function $T(r,g)$ to the both sides of \eqref{33}, then for any distinct value $a_{1},a_{2},\ldots,a_{l}\in\mathbb{C}\cup\{\infty\}$,
	\begin{align*}
		l-2\leq\sum_{i=1}^{l}\frac{\hat{N}(r,g=a_{i})}{T(r,g)}+\frac{o(T(r,g))}{T(r,g)}.
	\end{align*}
	That is,
	\begin{align}\label{3.1}
		\sum_{i=1}^{l}\left(1-\frac{\hat{N}(r,g=a_{i})}{T(r,g)}\right)\leq 2+\frac{o(T(r,g))}{T(r,g)}.
	\end{align}
	Hence, take liminf on the both sides of \eqref{3.1}, then
	\begin{align*}
		\sum_{i=1}^{l}(\delta(a,g)+\theta_{q,c}(a,g))\leq\sum_{i=1}^{l}\Theta_{q,c}(a,g)\leq2.
	\end{align*}
\end{proof}		
The value $a\in\mathbb{C}\cup\{\infty\}$ is called Hahn-Nevanlinna deficient value if $\Theta_{q,c}(r,g)>0$. The next result can be deduced from the Hahn difference defect relation by using similar proof method in [18, P.44], its proof is omitted here. \\

\begin{theorem}Suppose that $g$ is a non-constant meromorphic function with zero order. Then $g$ has at most countable number of Hahn-Nevanlinna deficient values.
\end{theorem}
%\subsection{Picard theorem for Hahn difference operators}
\indent For any $a\in\mathbb{C}\cup\{\infty\}$, if $\hat{n}_{q,c}(a,g)=O(1)$ (note that it is equivalent $\hat{N}_{q,c}(a,g)=O(\log r)$), the value $a$ is called a Hahn-Picard exceptional value of $g$. We note that $a$ is called a Hahn-Picard exceptional value of $g$ means that except for at most finite many points, the number of zeros of $g-a$ with multiplicities is not larger than the number of zeros of $\mathcal{D}_{q,c}g=0$ with multiplicities. We also remark that $\Theta_{q,c}(a,g)=1$ if transcendental function $g$ with a Hahn-Picard exceptional value $a$. Therefore, we will establish the analogue of Nevanlinna Picard theorem for Hahn difference operators base on the Theorem \ref{thm2} as below.\\		
\begin{theorem}
	Suppose that $g$ is a meromorphic function with zero order, and $0<|q|<1$. Then $g$ takes every complex value infinitely often, with at most two possible Hahn-Picard exceptional values.
\end{theorem}		
\begin{proof}
	Suppose that $g$ has three distinct Hahn-Picard exceptional values $a_{1},a_{2},a_{3}$, we aim for a contradiction. By the definition of Hahn-Picard exceptional values, then we have
	$g$ is not a non-constant polynomial and $\hat{N}_{q,c}(r,g=a_{i})=o(T(r,g)),i=1,2,3$. Applying the Theorem \ref{thm2},
	\begin{align*}
		T(r,g)\leq\sum_{i=1}^{3}\hat{N}_{q,c}(r,g=a_{i})=o(T(r,g))
	\end{align*}
	holds for all $r\in E_{1}$, where the set $E_{1}$ is of logarithmic density 1, which is impossible.
\end{proof}
%\subsection{Five-value theorem for Hahn difference operators}
\indent It is well-known that the celebrated Five-Value Theorem, obtained by R.~Nevanlinna in 1929, states that two non-constant meromorphic functions $g$ and $h$ must be identical if they share five distinct values $a_{1}, a_{2}, \dots, a_{5} \in \mathbb{C}\cup\{\infty\}$ ignoring multiplicity. This work advanced the development of the uniqueness theory for meromorphic functions. Next, we aim to establish an analogue of Five-value theorem in the Hahn sense. Before that, we need to explain what is the meaning of two functions Hahn-share value $a$.\\		
\indent For two meromorphic functions $g, h$ are of zero-order and $a\in\mathbb{C}\cup\{\infty\}$. Denote by $E_{g}(a)$ the set of zeros of $g=a$, which is a subset of $\mathbb{C}$. If $E_{g}(a)=E_{h}(a)$ except perhaps on the subset of $\mathbb{C}$ such that $$\hat{N}_{q,c}(r,g=a)-\hat{N}_{q,c}(r,h=a)=o(T(r,g)+T(r,h)),$$ then we say that $g,h$ Hahn-share the value $a$.
\begin{theorem}
	Suppose that $g,h$ are two non-constant meromorphic functions with zero order. If $g,h$ Hahn-share five distinct values $a_{1},\ldots,a_{5}\in\mathbb{C}\cup\{\infty\}$, then $g\equiv h$.
\end{theorem}		
\begin{proof}
	Assume that $g\not\equiv h$, we aim to get a contradiction. Applying the second main theorem for Hahn difference operators,
	\begin{align*}
		3(T(r,g)+T(r,h))&\leq\sum_{i=1}^{5}(\hat{N}_{q,c}(r,g=a_{i})+\hat{N}_{q,c}(r,h=a_{i}))\\
		&+o(T(r,g)+T(r,h))
	\end{align*}
	holds for all $r\in E_{1}$, where the set $E_{1}$ is of logarithmic density 1.			
	Since $g,h$ Hahn-share $a_{1},\ldots,a_{5}$,
	\begin{align*}
		\hat{N}_{q,c}(r,g=a_{i})-\hat{N}_{q,c}(r,h=a_{i})=o(T(r,g)+T(r,h)).
	\end{align*}			
	Hence,
	\begin{align}\label{t5.1}
		T(r,g)+T(r,h)\leq\frac{2}{3}\sum_{i=1}^{5}\hat{N}_{q,c}(r,g=a_{i})+o(T(r,g)+T(r,h))
	\end{align}
	holds for all $r\in E_{1}$.	At the same time, we can deduce from $g,h$ Hahn-share $a_{1},\ldots,a_{5}$ that
	\begin{align}\label{t5.2}
		\sum_{i=1}^{5}\hat{N}_{q,c}(r,g=a_{i})&\leq\hat{N}_{q,c}(r,g-h=0)\\
		&\leq T\left(r,\frac{1}{g-h}\right)\nonumber\\
		&\leq T(r,g)+T(r,h)+O(1), r\in E_{1}.
	\end{align}			
	Combining \eqref{t5.1} and \eqref{t5.2}, then
	\begin{align*}
		\frac{1}{3}(T(r,g)+T(r,h))\leq o(T(r,g)+T(r,h)),  r\in E_{1},
	\end{align*}
	which is a contradiction. So $g\equiv h$.
\end{proof}		

\section{Further results}
In 1927, Montel \cite{Montel} replaced \(x\) and \(y\) in Fermat's Diophantine equation \(x^{m}+y^{m}=1\) with functions \(f\) and \(g\), i.e. the functional equation
\begin{align}\label{4.1}
	f^{m}(z)+g^{m}(z)=1,
\end{align}
and proved that for \(m \geq 3\), all entire function solutions \(f\) and \(g\) of the Fermat equation \eqref{4.1} must be constant. The same conclusion can also be found in [23, Lemma 1]. Baker \cite{Baker} and Gross \cite{Gross} independently extended Montel's result, demonstrating that for \(m \geq 4\), equation \eqref{4.1} admits no non-constant meromorphic solutions. They further characterized the non-constant meromorphic solutions of \eqref{4.1} for the cases \(m = 2\) and \(m = 3\), respectively. In 1970, Yang \cite{Yang} considered the more general Fermat-type functional equation
\begin{align}\label{4.2}
	f^{n} + g^{m} = 1,
\end{align}
and proved that when \(\frac{1}{n} + \frac{1}{m} < 1\), the equation \eqref{4.2} admits no non-constant entire solutions. For equation \eqref{4.2}, many scholars have considered the scenario where \(n = m = 2\) and where there exists a specific functional relationship between \(f\) and \(g\), as seen in \cite{Li, LC, LQ, LY, TL, YL, YS}. The case of \(n = m = 3\) has also garnered considerable scholarly interest. In 1989, Yanagihara \cite{Yana} investigated the existence of meromorphic solutions to the Fermat difference equation
\begin{align}\label{4.3}
	f^{3}+f^{3}(z+c)=1
\end{align}
and demonstrated that the equation \eqref{4.3} admits no non-constant meromorphic solutions of finite order, which was also obtained by L\"{u} et al. in \cite{LV} by making use of the difference analogue of the logarithmic derivative lemma of finite order meromorphic functions.\\
\indent Due to the relationship between the Hahn difference operators and differential operators as well as difference operators, we consider the existence of meromorphic solutions with finite order for the equation
\begin{align}\label{4.4}
	f^{3}+(\mathcal{D}_{q,c}f(z))^{3}=1,~q\not\in\{0,1\},~c\in\mathbb{C},
\end{align}
and obtain the following Theorem \ref{thm6.5}. In order to prove it, we need a lemma as below, which can be seen in \cite{EF}.
\begin{lemma}\cite{EF}\label{l4.1}
	Let \(f\) be meromorphic and \(h\) be entire in \(\mathbb{C}\). Then, the following statements hold.\\
	\(\mathrm{(i)}\) If \(0<\rho(f),\rho(h)<\infty\), then \(\rho(f\circ h)=\infty\);\\
	\(\mathrm{(ii)}\) If \(\rho(f\circ h)<\infty\) and \(h\) is transcendental, then \(\rho(f)=0\).
\end{lemma}
\indent We consider equation \eqref{4.4} and obtain the following conclusion.
\begin{theorem}\label{thm6.5}
	The functional equation \eqref{4.4} does not admit non-constant meromorphic solutions of finite order.	
\end{theorem}
\begin{proof}
	Suppose that \(f(z)\) is a non-constant meromorphic solution of the equation \eqref{4.4}. Then, according to [7, Proposition 1],
	\begin{align}\label{pf4.1}
		f(z)=\frac{1}{2}\cdotp\frac{\{1+\frac{\wp'(h(z))}{\sqrt{3}}\}}{\wp(h(z))},~~~ \mathcal{D}_{q,c}f(z)=\frac{\eta}{2}\cdotp\frac{\{1-\frac{\wp'(h(z))}{\sqrt{3}}\}}{\wp(h(z))},
	\end{align}
	where \(\eta^{3}=1\), \(h\) is a non-constant entire function and \(\wp\) is the Weierstrass \(\wp\)-function satisfying
	\begin{align}\label{4}
		(\wp')^{2}=4\wp^{3}-1.
	\end{align}
	It can be obtained through a simple calculation that
	\begin{align}\label{pf4.2}
		&\{1+\frac{\wp'(h(qz+c))}{\sqrt{3}}\}{\wp(h(z))}\nonumber\\
		&=\{\eta[(q-1)z+c](1-\frac{\wp'(h(z))}{\sqrt{3}})+1+\frac{\wp'(h(z))}{\sqrt{3}}\}\wp(h(qz+c)).
	\end{align}
	\indent Combining \eqref{4} with \eqref{pf4.1}, then
	\begin{align*}
		\wp^{3}(h(z))=3f^{2}(z)\wp^{2}(h(z))-3f(z)\wp(h(z))	+1.
	\end{align*}
	\indent We suppose that \(\rho(f)<\infty\), then we can get \(\rho(\wp(h))<\infty\) because \(\rho(\wp)=2\). Furthermore, we can deduce that \(h\) is a polynomial by combining Lemma \ref{l4.1} and \(\rho(\wp)=2\).\\
	\indent By \eqref{4}, \((\wp')^{2}(z_{0})=-1\) if \(\wp(z_{0})=0\). Let \(\{z_{j}\}_{j=1}^{\infty}\) denote the set of all zeros of \(\wp\), where \(|z_{j}|\to\infty\) as \(j\to\infty\), and assume that \(h(a_{j,l})=z_{j}, l=1,2,...,\deg(h)\). Hence, we have \((\wp')^{2}(z_{j})=(\wp')^{2}(h(a_{j,l}))=-1\).\\
	\indent If \(\{a'_{j,l}\}_{j=1}^{\infty}\) is a infinite sub-sequence of \(\{a_{j,l}\}_{j=1}^{\infty}\) with respect to \(j\) such that \(\wp(h(qa'_{j,l}+c))=0)\), so \((\wp')^{2}(h(qa'_{j,l}+c))=-1\).\\
	\indent Differentiating equation \eqref{pf4.2} and substituting the values yields,
	\begin{align}\label{pf4.3}
		&\{1+\frac{\wp'(h(qa'_{j,l}+c))}{\sqrt{3}}\}\wp'(h(a'_{j,l}))h'(a'_{j,l})\nonumber\\
		&=[\eta(q-1)(1-\frac{\wp'(h(a'_{j,l}))}{\sqrt{3}})]\wp(h(qa'_{j,l}+c))\\
		&+\{[\eta(q-1)a'_{j,l}+c](1-\frac{\wp'(h(a'_{j,l}))}{\sqrt{3}})+(1+\frac{\wp'(h(a'_{j,l}))}{\sqrt{3}})\}\nonumber\\
		&\cdot\wp'(h(qa'_{j,l}+c))h'(qa'_{j,l}+c).\nonumber
	\end{align}
	\indent It is worth noting that equation \eqref{pf4.3} holds if and only if one of the following cases occurs
	\begin{align*}
		&\{1+i\frac{\sqrt{3}}{3}\}h'(a'_{j,l})=qB\{1-i\frac{\sqrt{3}}{3}\}h'(qa'_{j,l}+c),\\
		&\{1-i\frac{\sqrt{3}}{3}\}h'(a'_{j,l})=qB\{1+i\frac{\sqrt{3}}{3}\}h'(qa'_{j,l}+c),\\
		&h'(a'_{j,l})=-qBh'(qa'_{j,l}+c),	
	\end{align*}	
	where \(B=\eta[(q-1)a'_{j,l}+c]+1\). Since \(h\) and \(h(z+c)\) are polynomials of the same degree and there are infinitely many points \(a'_{j,l}\) satisfying \(|a_{j,l}|\to\infty\) as \(j\to\infty\), then,
	\begin{align*}
		&\{1+i\frac{\sqrt{3}}{3}\}h'(z)=qB\{1-i\frac{\sqrt{3}}{3}\}h'(qz+c),\\
		&\{1-i\frac{\sqrt{3}}{3}\}h'(z)=qB\{1+i\frac{\sqrt{3}}{3}\}h'(qz+c),\\
		&h'(z)=-qBh'(qz+c),	
	\end{align*}		
	which is impossible because of the value of \(\eta, B\). Therefore, \(\wp(h(qa'_{j,l}+c))=0\) holds for at most finitely many \(a'_{j,l}\)'s.	 Without loss of generality, there exists a sufficiently large positive integer \(J\), such that \(\wp(h(qa'_{j,l}+c))\neq0\) when \(j>J\) for \(l=1,2,...,\deg(h)\). It is evident from \eqref{pf4.2} that \(\wp(h(qa'_{j,l}+c))=\infty\) because \(\wp(h(a'_{j,l}))=0\) and \(\wp^{2}(h(a'_{j,l}))=-1\) for \(j>J\). Hence, combining this with \(O(\log r)=S(r,\wp(h))\), we have
	\begin{align}\label{pf4.4}
		N\left(r,\frac{1}{\wp(h(z))}\right)&\leq\overline{N}\left(r,\frac{1}{\wp(h(z))}\right)+2N\left(r,\frac{1}{h'(z)}\right)\nonumber\\
		&\leq\overline{N}(r,\wp(h(qz+c)))+2T(r,h')+O(\log r)\\
		&\leq\overline{N}(r,\wp(h(qz+c)))+S(r,\wp(h)).\nonumber
	\end{align}
	\indent According to \eqref{pf4.1} and \(\rho(\wp(h))<\infty\), then we have \(\rho(f)=\rho(\wp(h))\) and \(S(r,f)=S(r,\wp(h))\).	\\
	\indent On the other hand, it becomes evident that all zeros of \(f-1, f-\eta, f-\eta^{2}\) are of multiplicities at least 3 by expressing \((\mathcal{D}_{q,c}f)^{3}=f^{3}-1=(f-1)(f-\eta)(f-\eta^{2})\), where \(\eta\neq1\). Hence, application of the Second and the First main Theorems yields the following first and third inequalities, respectively,
	\begin{align*}
		2T(r,f)&\leq\sum_{m=1}^{3}\overline{N}\left(r,\frac{1}{f-\eta^{m}}\right)+\overline{N}(r,f)+S(r,f)\\
		&\leq\frac{1}{3}\sum_{m=1}^{3}N\left(r,\frac{1}{f-\eta^{m}}\right)+N(r,f)+S(r,f)\\
		&\leq2T(r,f)+S(r,\wp(h)),
	\end{align*}
	which implies that \(T(r,f)=N(r,f)+S(r,\wp(h))\) and \(m(r,f)=S(r,\wp(h))\).\\
	\indent According to the lemma of logarithmic derivative, then
	\begin{align}\label{pf4.5}
		m\left(r,\frac{1}{\wp(h)}\right)&=m\left(r,\frac{1}{2}\frac{1}{\wp(h)}\right)+O(1)\nonumber\\
		&\leq m(r,f)+m\left(r,\frac{\sqrt{3}}{6}\frac{\wp'(h)h'}{\wp(h)}\right)+T(r,h')+O(1)\\
		&=S(r,\wp'(h)).\nonumber
	\end{align}
	\indent It is worth noticing that all poles of \(\wp(h)\) have multiplicity \(2t~(t\geq1)\) because each pole of \(\wp\) is of multiplicity 2. Combining \eqref{pf4.4} and \eqref{pf4.5}, then
	\begin{align*}
		T(r,\wp(h))&=T\left(r,\frac{1}{\wp(h)}\right)+O(1)\\
		&=m\left(r,\frac{1}{\wp(h)}\right)+N\left(r,\frac{1}{\wp(h)}\right)+O(1)\\
		&\leq\overline{N}(r,\wp(h(qz+c)))+S(r,\wp(h))\\
		&\leq\frac{1}{2}N(r,\wp(h(qz+c)))+S(r,\wp(h))\\
		&\leq\frac{1}{2}T(r,\wp(h(qz+c)))+S(r,\wp(h))\\
		&=\frac{1}{2}T(r,\wp(h(z)))+S(r,\wp(h))+O(r^{\rho(\wp(h))-1+\varepsilon}),
	\end{align*}
	where the last estimation can be found in [8, Theorem 2.1] with \(\varepsilon>0\) sufficiently small. This is a contradiction.
\end{proof}

\section{Question}
\indent The Wiman-Valiron theory is often employed to characterize the upper bounds of growth for solutions of complex linear differential equations. Besides, the Wiman-Valiron theory for difference operators and \(q\)-difference operators were established by Ishizaki-Yanagihara \cite{IY} and Wen-Ye \cite{WY}, respectively.\\
\indent  Therefore, we pose the following question: Whether a Hahn difference version of the Wiman-Valiron theory can be established?

\section*{Author Contributions} All authors wrote and completed the main part of this article and corrected the main theorems. All authors gave final approval for publication.

\section*{Statements and Declarations}
\(\bullet\) This research work was supported by the National Natural Science Foundation of China (Grant No. 12261023), the Graduate Research Fund Project of Guizhou Province (2024YJSKYJJ186), and the Guizhou Key Laboratory of Advanced Computing(Grant No. QianKeHe Platform ZSYS(2025)004).\\
\(\bullet\) No conflict of interest.\\
\(\bullet\) Compliance with ethical standards.

\end{document}